\documentclass[reqno]{amsart}[a4paper]
\usepackage{latexsym,amsfonts,amssymb,epsfig}
\usepackage{hyperref} %% autoreferences, works only with PdfLaTeX
\usepackage{color}

\newcommand{\N}{\mathbb N}            % natural numbers
\newcommand{\Z}{\mathbb Z}            % integers
\newcommand{\R}{\mathbb R}
\newtheorem{Th}{Theorem}
\newtheorem{Corr}[Th]{Corollary}%[section]
\theoremstyle{Conjecture}
\newtheorem{Conjecture}{Conjecture}
\renewcommand{\AA}{\mathbf A}         % abelian   variety
\newcommand{\NN}{\mathbf N}           % nilpotent variety
\DeclareMathOperator{\rank}{rank}
\DeclareMathOperator{\Exp}{Exp}
\def\a{\alpha}
\def\g{\gamma}

\begin{document}
\title{Schreier's type formulae and two scales for growth of Lie algebras and groups}
\author{Victor Petrogradsky}
\address{Department of Mathematics, University of Brasilia, 70910-900 Brasilia DF, Brazil}
\email{petrogradsky@rambler.ru}
%\thanks{The author was partially supported by grants CNPq~309542/2016-2, DPI-UnB~04/2019}
%University of Brasilia, 70910-900 Brasilia DF, Brazil, petrogradsky@rambler.ru, pitersky@unb.br\\
%ORCID 0000-0002-0182-5149}
%\\ \vspace{-0.3cm} %\center{ \bf\today}
\dedicatory{Dedicated to 70th anniversary of Vesselin Drensky}
%\\ ORCID 0000-0002-0182-5149}

\begin{abstract}
Let $G$ be a free group of rank $n$ and $H\subset G$ its subgroup of finite index.
Then $H$ is also a free group and the rank $m$ of $H$ is determined by Schreier's formula
$m-1=(n-1)\cdot|G:H|.$

Any subalgebra of a free Lie algebra is also free.
But a straightforward analogue of Schreier's formula for free Lie
algebras does not exist, because any subalgebra of finite codimension has an infinite number of generators.

But the appropriate Schreier's formula for free Lie algebras exists in terms of formal power series.
There exists also a version in terms of exponential generating functions.
This is a survey on how these formulas are applied to study
1) growth of finitely generated Lie algebras and groups
and 2) the codimension growth of varieties of Lie algebras.
First, these formulae allow to specify explicit formulas for generating functions
of respective types for free solvable (or more generally, polynilpotent) Lie algebras.
Second, these explicit formulas for generating functions are used to derive asymptotic for these two types of the growth.
These results can be viewed as analogues of the Witt formula for free Lie algebras and groups.
In case of Lie algebras, we obtain two scales for respective types of growth.
We also shortly mention the situation on growth for other types of linear algebras.
\end{abstract}

\maketitle

%*********************************************************************************************
\section {Analogue of Schreier's formula for free Lie algebras}
Denote the ground field by $K$. Let $X$ be an at most countable set supplied with a
{\it weight function} $\mathop{\mathrm {wt}}\nolimits: X\to\mathbb N$, namely we assume that
$$\displaystyle X=\mathop{\cup}_{i=1}^\infty X_i;\quad X_i=\{x\in X\
|\mathop{\mathrm {wt}}\nolimits x=i\},\quad |X_i|<\infty,\ i\in\mathbb N.$$
We say that $X$ is {\em finitely graded}.
Assume that an algebra $A$ is generated by $X$, then we naturally define the weight of a monomial $a\in A$.
Let $Y$ be a set of monomials in $X$, then one defines the {\it Hilbert-Poincar\'e series} of $Y$
(with respect to $X$), see e.g.~\cite{Drensky}
$${\mathcal H}_X(Y)={\mathcal H}_X(Y,t):=\mathop{\sum}_{i=1}^\infty |Y_i|t^i;\qquad
Y_i=\{y\in Y\ |\mathop{\mathrm {wt}}\nolimits y=i\},\ i\in\mathbb N.$$
Consider a subspace $V\subset A$, then we define ${\mathcal H}_X(V)$ using a homogeneous basis of $\mathop{\rm gr}\nolimits V$;
where $\mathop{\rm gr}\nolimits V$ being the associated graded space.

Next, we introduce the operator ${\mathcal E}$ on power series
$\phi(t)=\sum_{n=1}^\infty b_n t^n$, $b_n\in\{0,1,2,\dots\}$ (see~\cite{Pe99int,Pe00archiv}):
$${\mathcal E}\,:\,\phi(t)=\sum_{n=1}^\infty b_n t^n\ \longrightarrow\
   {\mathcal E}(\phi(t)):=\sum_{n=0}^\infty a_n t^n
   =\prod_{n=1}^\infty \frac 1{(1-t^n)^{b_n}}. $$
Assume that $L$ is a Lie algebra generated by $X$ and  $U(L)$ its universal enveloping algebra. Let
$$
{\mathcal H}_X(L,t)=\sum_{n=1}^\infty b_n t^n,\quad
{\mathcal H}_X(U(L),t)=\sum_{n=0}^\infty a_n t^n.
$$
One has a well-known formula that explains importance of the operator above
${\mathcal H}_X(U(L))={\mathcal E}({\mathcal H}_X(L))$~\cite{Ufn}.
The following is the natural analogue of Schreier's formula, introduced by the author.
For basic facts on free Lie (super)algebras see~\cite{Ba,BMPZ}.

\begin{Th}[\cite{Pe00archiv}]\label{T1}
Assume that $L$ is a free Lie algebra generated by a finitely graded set $X$.
Let $H$ be a subalgebra and $Y$ is a set of its free generators. Then
$${\mathcal H}(Y,t)-1=({\mathcal H}(X,t)-1)\cdot{\mathcal E}({\mathcal H}(L/H),t).$$
\end{Th}

Let $L$ be a Lie algebra. One defines the {\it lower central series} as
$L^1=L$, $L^{i+1}=[L,L^i]$ for  $i=1,2,\dots$.
Now, $L$ is {\em nilpotent} of class $s$ iff $L^{s+1}=\{0\}$ while $L^s\not=\{0\}$.
All Lie algebras nilpotent of class at most $s$ form the variety denoted by $\NN_s$.
A Lie algebra $L$ is called {\em polynilpotent} with a tuple of integers $(s_q,\dots,s_2,s_1)$
iff there exists  a chain of ideals
$$\{0\}=L_{q+1}\subset L_q\subset\dots\subset L_2\subset L_1=L,\qquad L_n/L_{n+1}\in \NN_{s_n}, \ n=1,\ldots,q.$$
All polynilpotent Lie algebras with a fixed tuple form
a variety denoted by $\NN_{s_q}\dots \NN_{s_2}\NN_{s_1}$.
In the particular case $s_q=\cdots=s_1=1$, one obtains the variety ${\AA}^q$,
consisting  of {\em solvable} Lie algebras of length at most $q$.
On the other hand, a polynilpotent Lie algebra is solvable.
Moreover, the free polynilpotent Lie algebras, a tuple being fixed,
provide interesting examples of solvable Lie algebras.
The definitions in case of group theory are similar.
Let $G$ be a group, denote by $\{\gamma_n(G)|n=1,2,\ldots\}$,
terms of the {\it lower central series} (warning: below $\gamma$ has also a different meaning!).

Suppose that $L$ is the free Lie algebra of rank $k$,
$L=\oplus_{n=1}^\infty L_n$ its natural grading,
and $G$ the free group of rank $k$. Then the lower central series factors $\g_n(G)/\g_{n+1}(G)$
are free abelian groups and their ranks are given by the classical Witt formula~\cite{Ba}:
\begin{equation}\label{witt}
\psi_k(n):=\rank_\Z \g_n(G)/\g_{n+1}(G)=\dim_K L_n=
\frac 1n\sum_{a|n}k^a\mu\left(\frac na\right)\approx\frac{k^n}n, \quad n\in \N,
\end{equation}
where $\mu(*)$ is the M\"obius function.
Theorem~\ref{T1} allows to derive the following explicit formulas.
\begin{Th}[\cite{Pe00archiv,Pe02}]\label{TX}
Consider the free polynilpotent Lie algebra $L=F({\bf{N}}_{s_q}\cdots{\bf{N}}_{s_1},k)$ of finite rank $k\ge 2$, where $q\ge 1$.
Set $\beta_0(z):=0$, $\alpha_0(z):=kz$, and define the following functions  recursively
  \begin{align*}
  \begin{split}
  \beta_i(z)&:=\beta_{i-1}(z)+\sum_{m=1}^{s_i}\frac 1m
           \sum_{a|m}\mu\left(\frac ma\right)
           \left(\alpha_{i-1}(z^{m/a})\right)^a\!\!,\\
  \alpha_i(z)&:=1+(kz-1)\cdot{\mathcal E}(\beta_i(z)),
  \end{split} \qquad 1\le i\le q.
  \end{align*}
Then ${\mathcal H}(L,z)=\beta_q(z).$
\end{Th}
For example, we have a particular case.
\begin{Corr}[\cite{Pe00archiv,Pe02}]
Let $L:=F({\bf{A}}{\bf{N}}_d,k)$. %% Denote $\psi_k(m):=\frac 1m \sum_{a|m} \mu(m/a)k^a$.
Then
$$
{\mathcal H}(L,z)=\psi_k(1)z+\dots+\psi_k(d)z^d+1+
       \frac{kz-1}
       {(1-z)^{\psi_k(1)}%(1-z^2)^{\psi_k(2)}
       \cdots(1-z^d)^{\psi_k(d)}}.
$$
\end{Corr}

%****************************************************************************************************
\section {Scale 1 and Growth of free solvable (polynilpotent) finitely generated Lie algebras and groups}

Now we describe applications of the analogue of Schreier's formula for free Lie algebras
(Theorem~\ref{T1} and its application Theorem~\ref{TX})
to specify the growth of free solvable (more generally, polynilpotent) Lie algebras and groups of finite rank.

Assume that $L$ is a relatively free algebra of some multihomogeneous variety of
(associative, or Lie) algebras,
generated by $X=\{x_1,\dots,x_k\}$.
Then we have a natural grading $L=\mathop{\oplus}_{n=1}^\infty L_n$ by degree in~$X$.
One defines the {\em growth function} with respect to $X$ as $\gamma_L(X,n):=\sum_{s=1}^n \dim_K L_s$.

\begin{Th}[Berele, \cite{Ber82}]\label{Tberele}
The growth function of a finitely generated associative PI-algebra is bounded by a polynomial function.
\end{Th}

For more details on proofs of this important result see~\cite{Drensky,KraLen}.
But the growth of finitely generated Lie PI-algebras is more complicated.
In this case, the author introduced scale~\eqref{scale1} of functions of intermediate growth and
suggested that it is complete in the sense of Conjecture~1 below.
Define functions
\begin{align*}
\begin{split}
\ln^{(0)}x:=x,\qquad  &\ln^{(s+1)}x:=\ln(\ln^{(s)}x),\\
\exp^{(0)}x:=x,\qquad &\exp^{(s+1)}x:=\exp(\exp^{(s)}x),
\end{split} \qquad\qquad
s=0,1,2,\dots
\end{align*}
Recall standard notations
$f(x)\approx g(x)$, $x\to\infty$,
denotes that $\mathop{{\rm lim}\,}_{x\to\infty}f(x)/g(x)=1$;
$f(x)=o(g(x))$ means that beginning with some number
$f(x)=\alpha(x)g(x)$ and $\mathop{{\rm lim}\,}_{x\to\infty}\alpha(x)=0$.
%also $f(x)=O(g(x))$ iff $f(x)=\beta(x)g(x)$ where $\beta(x)$ is bounded.
Also, $\zeta(x)$ is the zeta function.

Consider the {\bf scale 1} consisting of a
a series of functions $\Phi^q_\alpha (n)$, $q=1,2,3,\dots$
of a natural argument with a parameter $\alpha\in\R^+$:
\begin{align}\label{scale1}
\text{\textbf {scale 1:}}\qquad
\begin{split}
 \Phi^1_\a(n)&:=\a,\\
 \Phi^2_\a(n)&:=n^\a,\\
 \Phi^3_\a(n)&:=\exp(n^{\a/(\a+1)}),\\
 \Phi^q_\a(n)&:=\exp\bigg({n\over (\ln^{(q-3)}n)^{1/\a} }\bigg),\quad
        \quad q=4,5,\dots.
\end{split}
\end{align}

Now, we specify the growth of free solvable (more generally, polynilpotent)
Lie algebras of finite rank with the respect to scale~\eqref{scale1} by giving the following asymptotic.

\begin{Th}[\cite{Pe99int}]\label{T4}
Consider the free polynilpotent Lie algebra $L:=F({\bf{N}}_{s_q}\cdots{\bf{N}}_{s_1},k)$ of finite rank $k\ge 2$, where $q\ge 2$,
generated by $X=\{x_1,\dots,x_k\}$. Then
$$\displaystyle
\gamma_{L}(X,n)=
  \left\{
  \begin{array}{ll}
     \displaystyle\frac {A+o(1)}{N!}\,n^N,
       & q=2, \\[10pt]
     \displaystyle\exp\Big((C{+}o(1))\,n^{{N}/(N+1)}\Big),
       & q=3,  \\[10pt]
\displaystyle\exp\Biggl( (B^{1/N}{+}o(1))\,{ \frac{n}{(\ln^{(q-3)}n)^{1/N}}}
\Biggr),\quad
       & q\ge 4,
  \end{array}
  \right.\qquad n\to\infty;
$$
where the constants are:
\begin{eqnarray*}
&&N:=s_2\dim_K F({\bf{N}}_{s_1},k),\quad
A:=\frac 1{s_2} \left(\frac{k-1}{\prod_{q=2}^{s_1}
q^{\psi_k(q)}} \right)^{s_2},\\
&&B:=s_3A\zeta(N{+}1),\qquad\qquad
C:=(1+1/N)(BN)^\frac 1{1+N}.
\end{eqnarray*}
\end{Th}

This result has the following application to group theory. Let $G$ be a group.
Due to Lazard, one constructs the related Lie algebra~\cite{Laz54}:
\begin{align*}
L_K(G):=\mathop{\oplus}\limits_{n=1}^\infty (\gamma_n(G)/\gamma_{n+1}(G))\otimes_{\Z}K.
\end{align*}
If $G$ is a free polynilpotent group, then $L_K(G)$ is the free polynilpotent
Lie algebra of the same rank and with the same tuple, moreover,
the lower central series factors are free abelian groups (A.~Shmelkin~\cite{Shmelkin}).

\begin{Corr}[\cite{Pe99int,Pe02}]
\label{Ccentr}
Let $G=G(\NN_{s_q}\dots \NN_{s_1},k)$, $q\geq 2$, be the free polynilpotent group  of rank $k$.
Consider ranks of the lower central series factors, namely, denote
%Let $\gamma_n(G)$, $n=1,2,\dots $, be the terms of the lower central series. Denote
$b_n:=\rank_\Z \gamma_n(G)/\gamma_{n+1}(G)$, $n\ge 1$. Then we get the asymptotic:
$$ b_n=
  \left\{
  \begin{array}{ll}
     \displaystyle\frac {A+o(1)}{(N-1)!}\,n^{N-1},
       & q=2, \\[10pt]
     \displaystyle\exp\Big((C{+}o(1))\,n^{{N}/(N+1)}\Big),
       & q=3,  \\[10pt]
     \displaystyle\exp\Biggl( (B^{1/N}{+}o(1))\,{ n\over (\ln^{(q-3)}n)^{1/N}} \Biggr),\quad
       & q\ge 4,
  \end{array}
  \right.\qquad n\to\infty;
$$
where $N,A,B,C$ are the same as in Theorem~\ref{T4}.
\end{Corr}

Observe that just by setting $s_q=\cdots=s_1=1$,
Theorem~\ref{T4} and its Corollary~\ref{Ccentr} are turned into results on the free solvable
Lie algebra and group of rank $k$ and length $q$.
A similar observation is valid concerning the results  (e.g.~Theorem~\ref{TNN2}) on the codimension growth below.

M.I.Kargapolov raised problem 2.18  in~\cite{KouTet67}
to specify the lower central series ranks for free polynilpotent finitely generated  groups.
Exact recursive formulae  were  given by Egorychev~\cite{Egor}.
We suggest another answer to this problem by specifying the asymptotic behaviour of that ranks.
We consider to view the asymptotic of Theorem~\ref{T4} and its Corollary~\ref{Ccentr}
as an analogue of the Witt formula~\eqref{witt}, now for the free solvable (more generally, polynilpotent) Lie
algebras and groups.
\medskip

Since the growth of finitely generated solvable Lie algebras is intermediate~\cite{Licht84},
the respective generating function converges in the unit circle.
It is important to study an asymptotic of the generating function when
we approach the unit circumference from inside.
Since the coefficients of the series are real nonnegative numbers,
it is sufficient to study this behaviour while $z=t\to 1-0$.
The crucial step to prove Theorem~\ref{T4} is the following asymptotic of the generating function  ${\mathcal H}_X(L,t)$.
\begin{Th}[\cite{Pe99int}]\label{T5}
Let $L:=F({\bf{N}}_{s_q}\cdots{\bf{N}}_{s_1},k)$  be
the free polynilpotent Lie algebra of finite rank $k\ge 2$, where $q\ge 2$, generated by $X=\{x_1,\dots,x_k\}$. Then
$$
\begin{array}{ll}
\displaystyle\lim_{t\to 1{-}0\,}(1-t)^N {\mathcal H}_X(L,t)=A, & q=2, \\[8pt]
\displaystyle\lim_{t\to 1{-}0\,}(1-t)^N \ln^{(q-2)}({\mathcal
H}_X(L,t))=s_3\,\zeta(N{+}1)\,
A,\quad & q\ge 3;
\end{array}
$$
where $N,A$ are the same as in Theorem~\ref{T4}.
\end{Th}

Now, we describe an important idea to prove Theorem~\ref{T5}.
The explicit formula of the generating function (Theorem~\ref{TX})
shows that ${\mathcal H}_X(L,t)$ is "roughly speaking"
a $(q-1)$-iteration of ${\mathcal E}$ applied to $\beta_1(z)$.
We can easily describe the first application:
\begin{align*}
  \beta_1(z)&=\psi_k(1)z+\psi_k(2)z^2+\dots+\psi_k(s_1)z^{s_1};\\
%\quad \psi_k(m):=\frac {1}{m} \sum_{a|m} \mu(m/a)k^a;\\
  {\mathcal E}(\beta_1(z))&=
    \frac {1}
{(1-z)^{\psi_k(1)}%%(1-z^2)^{\psi_k(2)}
\cdots(1-z^{s_1})^{\psi_k(s_1)}}
    \approx\frac{\mu} {(1-t)^M},\quad \text{as}\quad  t\to 1{-}0;\\
 \text{where}\quad   M&=\psi_k(1)+\dots+\psi_k(s_1)=\dim_K F({\bf{N}}_{s_1},k),\qquad
           \mu=\prod_{q=2}^{s_1} q^{-\psi_k(q)}.
\end{align*}

Next, we use the fact that ${\mathcal E}$ is "approximately" the exponent, thus ${\mathcal H}_X(L,t)$ behaves like
$$ \exp^{(q-2)}\left(\frac{\bar\mu+o(1)}{(1-t)^N} \right),\qquad t\to 1{-}0\,. $$
Another idea in the proof of Theorem~\ref{T4} is to specify a connection between the growth
of a function analytic in the unit circle with asymptotic of its coefficients~\cite{Pe99int}.

A similar version of Schreier's formula for
free Lie {\it super}algebras was established as well~\cite{Pe00archiv,Pe02}.
Also, the asymptotic of Theorem~\ref{T4} was extended to the case of free solvable (polynilpotent) Lie
{\it super}algebras of finite rank~\cite{KlPe05}.

Below we see that the scale~\eqref{scale2} for the superexponential codimension growth for varieties of Lie algebras
is rather complete (Theorem~\ref{Tscale2}).
So, we conjecture that the scale~\eqref{scale1} for the intermediate growth of finitely generated Lie PI-algebras is complete as well.
\begin{Conjecture}[\cite{Pe96}]
Let $L$ be a finitely generated Lie PI-algebra. Then there exist numbers $q,N_0$ such that
$$
\gamma_L(n)\le \exp\bigg( \frac{n}{ \ln^{(q)}n } \bigg),\quad n\ge N_0.
$$
\end{Conjecture}
This bound was confirmed for almost solvable (more generally, almost polynilpotent)
Lie algebras~\cite{KlPe05usp,KlPe09contempmath}.

\section {Exponential analogue of Schreier's formula for free Lie algebras}
Let us consider {\it complexity functions}, refereed to also as {\it exponential generating functions} in combinatorics.
Assume that we are given a set $A$ of monomials in $X=\{x_i\,|\,i\in\mathbb N \}$.
Consider a set of distinct elements $\widetilde{X} =\{x_{i_1},\dots,x_{i_n}\}\subset X$,
denote by  $P_n(A,\widetilde{X})$ the set of all multilinear  elements of degree $n$ on $\widetilde{X}$ belonging to $A$.
Suppose that the number of these elements $c_n(A,\widetilde{X})$
does not depend on the choice of $\widetilde{X}$, but depends only on $n$.
In this case, we denote $c_n(A):=c_n(A,\widetilde{X})$ and say that
$A$ is $X$-{\em uniform} and define the {\it complexity function with respect to} $X$:
\begin{equation}\label{complexity}
{\mathcal C}_X(A,z):=\sum_{n=1}^\infty \frac{c_n(A)}{n!}z^n,\quad z\in \mathbb C.
\end{equation}
(the sum is taken from $n=0,\  c_0=1$ for associative algebras and groupoids with unity).
Remark that $A$ need not consist of multilinear elements.
We also omit the variable $z$ and (or) the set $X$  and write ${\mathcal C}_X(A,z)={\mathcal C}(A)$.
These definitions are naturally extended to
algebras, subspaces, and their dimensions  with respect to their  generating sets.
We illustrate this notion by examples. Let $X$ be a countable set and
$A=A(X)$, $L=L(X)$ the free associative and Lie algebras, respectively. Then
\begin{align*}
{\mathcal C}_X(A,z)&=\sum_{n=0}^\infty z^n=\frac 1{1-z},\\
{\mathcal C}_X(L,z)&= \sum_{n=1}^\infty \frac {z^n}n =-\ln(1-z),\\
{\mathcal C}_X({\bf{N}}_s,z)&= \sum_{n=1}^s \frac {z^n}n.
\end{align*}

The author established an exponential analogue of Schreier's formula for free Lie algebras as follows.
\begin{Th}[\cite{Pe99JMSciSch_exp}]  \label{Texp}
Assume that $L$ is the free Lie algebra generated by a countable generating set $X$.
Assume that $H$ is an $X$-uniform subalgebra.
Then $H$ has an $X$-uniform set of free generators $Y$ and
$${\mathcal C}_X(Y,z)-1=(z-1)\cdot\exp({\mathcal C}_X(L/H,z)).$$
\end{Th}

%***************************************************************************************
\section {Scale 2 for the codimension growth of Lie PI-algebras}

For the theory of varieties of associative and Lie algebras see~\cite{Ba,Drensky,GiaZai}.
Let ${\bf{V}}$ be a variety of Lie  algebras,
and $F({\bf{V}},X)$ its free algebra generated by $X=\{x_i | i\in\mathbb N\}$.
Let  $P_n({\bf{V}})\subset F({\bf{V}},X)$  be the subspace  of all
multilinear elements in $\{x_1,\dots,x_n\}$ and consider the {\em codimension growth sequence}
$c_n({\bf{V}})=c_n(F({\bf{V}},X),X):=\dim_KP_n({\bf{V}})$, $n=1,2,\dots$.

In case of associative algebras the fundamental fact is as follows:
\begin{Th}[Regev, \cite{Reg72}; Latyshev~\cite{Lat63}] \label{ThReg}
Let an associative algebra $A$ satisfies a nontrivial identical relation of  degree $d$. Then
$c_n(A)\leq C^n$, $n\ge 1$; where $C:=(d-1)^2.$
\end{Th}
Another crucial fact on the codimension growth of associative algebras in characteristic zero
is that the {\it exponent}, defined as
$\Exp(A):=\mathop{\lim}\limits_{n\to \infty} \sqrt[n]{c_n(A)}$
always exists and is integral~\cite{GiaZa99}.

Now we start discussing the codimension growth for Lie algebras.
The integrality of the exponent of the codimension growth
for finite dimensional Lie algebras over a field of characteristic zero was proved by Zaitsev~\cite{Zai02}.
In general, the exponents for Lie algebras are not always integral~\cite{ZaiMi99,GiaZa13}.
Moreover, the codimension growth in case of Lie algebras is more versatile.
Unlike the associative case, the codimension growth of a rather small variety $\AA\NN_2$
is overexponential (Volichenko \cite{Vol}).
On the other hand, the following upper bound was found.
\begin{Th}[Grishkov \cite{Grishkov88}]
Let $L$ be a Lie algebra satisfying a nontrivial identity.
Then for any $r>1$ there exists $N_0$ such that
$$ c_n(L)\le \frac {n!}{r^n},\quad n\ge N_0.$$
\end{Th}
Razmyslov introduced the {\em complexity functions}~\eqref{complexity} and  reformulated the upper bound as follows.%~\cite{Razmyslov}.
\begin{Th}[Razmyslov~\cite{Razmyslov}]
Let  ${\bf{V}}$  be a nontrivial variety of Lie algebras. Then the complexity function
${\mathcal C}({\bf{V}},z)$ is an entire function of complex variable.
\end{Th}
The author established a better and optimal general bound for the series that allowed to prove an upper bound for
the codimension growth sequence as well.
The estimate of the first item was recently obtained in~\cite{Pe21}.
\begin{Th}[\cite{Pe95,Pe97msb,Pe99isr,Pe21}]\label{Tscale2}
Let $L$ be a Lie algebra satisfying a nontrivial identity of degree $m\ge 4$. Then
\begin{enumerate}
\item The following coefficientwise bound for the series holds:
$$
\mathcal{C}(L,z)\prec z\underbrace{\exp(z \exp(\ldots(z\exp(z\exp}\limits_{m-2 \text{\rm{ times }} \exp } (z)))  \ldots)).
$$
\item there exists an infinitesimal such that
$$ c_n(L) \le \frac{n!} {(\ln^{(m-3)}n)^n}(1+o(1))^n, \quad n\to\infty. $$
\end{enumerate}
\end{Th}

Thus, we have a vast area of overexponential growths for Lie algebras,
lying between the exponent and the factorial functions.
To describe such a growth we introduce the {\bf scale 2} consisting of a series of functions
$\Psi^q_\a(n)$, $q=2,3,\dots$, with a real parameter $\a$~\cite{Pe95}:
\begin{equation}
\text{\textbf {scale 2:}}\qquad
\label{scale2}
\Psi^q_\a(n):=
\left\{
\begin{array}{lll}
\displaystyle
\big(n!\big)^\frac{\scriptstyle\a-1}{\scriptstyle\a},\quad & \a>1,\quad & q=2;\\
\displaystyle
\frac{n!}{(\ln^{(q-2)}n)^{n/\a}},\quad & \a>0,\quad & q=3,4\dots
\end{array}
\right.
\end{equation}

The upper bounds of Theorem~\ref{Tscale2} are "adequate" and the scale~\eqref{scale2}
for the codimension growth of Lie PI-algebras is complete,
because the free solvable Lie algebras
do have such an asymptotic behaviour, see Theorem~\ref{TNN2} below.
Actually we obtain a more fine scale formed by a series of
functions with two real parameters $\alpha,\beta$:
\begin{equation}
\label{scaleCo2}
\text{\textbf {scale 2$'$:}}\qquad
\Psi^q_{\alpha,\beta}(n):=
\left\{
\begin{array}{lll}
\displaystyle
\big(n!\big)^\frac{\scriptstyle\alpha-1}{\scriptstyle\alpha}\beta^{n/\alpha},
\quad & \alpha\ge 1,\ \beta>0;\quad & q=2;\\[7pt]
\displaystyle
\frac{n!\cdot (\beta/\alpha)^{n/\alpha}}{(\ln^{(q-2)}n)^{n/\alpha}},
\quad & \alpha>0,\ \beta>0; & q=3,4\dots
\end{array}
\right.
\end{equation}
Observe that in terms of scale~\eqref{scaleCo2} the exponential growth is a subcase of level $q=2$ when $\a=1$.

\subsection{The codimension growth for another classes of {\it linear algebras}}
The codimension growth of arbitrary linear algebras can be weird~\cite{GiaMiZai08}.
The varieties of absolutely free (commutative, or anticommutative) algebras
have Schreier's type formulae in terms of generating function of both types,
the regular generating functions and exponential generating functions (i.e. complexity functions),
i.e. we have natural analogues of both, Theorem~\ref{T1} and Theorem~\ref{Texp}), see~\cite{Pe05}.
We describe both, the generating functions and the growth functions, for two kinds of growth, for
different versions of nilpotency and solvability for three types of linear algebras above~\cite{Pe05}.
But we do not get an analogue of Theorem~\ref{Tproduct} and
these results do not lead us to something like scale~1 and scale~2 for the respective types of growth~\cite{Pe05}.

It was recently shown that the same scale~\eqref{scale2}  stratifies the ordinary codimension growth of
{\it Poisson PI-algebras},
but here it is essential to assume that a Poisson algebra satisfies a nontrivial Lie identical relation~\cite{Pe21}.
If a Poisson algebra is satisfying so called mixed identities only,
then the ordinary codimension growth has a factorial behaviour~\cite{Rats14,Pe21}.

The codimension growth for {\it Jordan PI-algebras} can be  overexponential \cite{Drensky87,GiaZe11}.
We conjecture that something like scale 2 (see~\eqref{scale2}) should appear in case of Jordan PI-algebras as well.

\section {Explicit formulae for complexity functions and asymptotic for Lie PI-algebras}

Let ${\bf{M}}, {\bf{V}}$ be varieties of Lie algebras.
Their {\it product}  ${\bf{M}}\cdot{\bf{V}}$ is the class of all Lie algebras $L$ such that there exists an ideal $H\subset L$
satisfying $H\in {\bf{M}}$ and $L/H\in {\bf{V}}$, see~\cite{Ba}.
Using the exponential Schreier's formula (Theorem~\ref{Texp}) the following explicit formula was proved.
\begin{Th}[\cite{Pe99JMSciSch_exp}] \label{Tproduct}
Let ${\bf{M}}\cdot{\bf{V}}$ be the product of varieties of
Lie algebras, where ${\bf{M}}$ is multihomogeneous. Then
$$ {\mathcal C}({\bf{M}}\cdot{\bf{V}},z)={\mathcal C}({\bf{V}},z)+{\mathcal C}({\bf{M}},1+(z-1)\exp({\mathcal C}({\bf{V}},z))). $$
\end{Th}

Roughly speaking, the formula says that ${\mathcal C}({\bf{M}}\cdot{\bf{V}},z)$ is "almost"
a composition of three functions ${\mathcal C}({\bf{M}})\circ\exp\circ\,{\mathcal C}({\bf{V}})$.
The variety ${\bf{V}}={\bf{N}}_{s_q}\cdots{\bf{N}}_{s_1}$ can be viewed as the product ${\bf{V}}={\bf{N}}_{s_q}\cdots{\bf{N}}_{s_2}\cdot{\bf{N}}_{s_1}$.
As application, the following explicit formula was derived.

\begin{Th}[\cite{Pe99JMSciSch_exp,Pe02}]\label{TcodimNN}
Consider the variety of polynilpotent Lie algebras
${\bf{V}}:={\bf{N}}_{s_q}\cdots{\bf{N}}_{s_1}$, $q\ge 1$. Define functions
  \begin{align*}
  \begin{split}
  \beta_1(z)&:=\displaystyle \sum_{m=1}^{s_1}\frac {z^m}m,\\
  \beta_i(z)&:=\displaystyle \beta_{i-1}(z)+\sum_{m=1}^{s_i}
         \frac{(1+(z-1)\exp(\beta_{i-1}(z)))^m}m, \qquad 2\le i\le q.
  \end{split}
  \end{align*}
Then ${\mathcal C}({\bf{V}},z)=\beta_q(z).$
\end{Th}
Consider particular cases.
\begin{Corr} Fix $d\in \N$. Then
$$ {\mathcal C}({\bf{A}}{\bf{N}}_d,z)=z+\frac{z^2}2+\dots+\frac{z^d}d+1+
       (z-1)\exp\left(z+\frac{z^2}2+\dots+\frac{z^d}d\right). $$
\end{Corr}
\begin{Corr}[\cite{Pe97msb}] Fix $q\in\N$.
Consider the variety of solvable Lie algebras of length $q$, denoted as ${\bf{A}}^q$.
Set $\beta_1(z)=z$, and $\beta_{i+1}(z)=\beta_i(z)+1+(z-1)\exp(\beta_i(z))$ for $i=1,2,\dots, q-1$.
Then $${\mathcal C}({\bf{A}}^q,z)=\beta_q(z).$$
\end{Corr}
\begin{Corr}[\cite{Pe97msb}] Fix $c\in \N$. Then
$$ {\mathcal C}({\bf{N}}_c{\bf{A}},z)=z+ \sum^c_{m=1}\frac 1m \Big(1+(z-1)\exp(z)\Big)^m. $$
\end{Corr}
Complexity functions are useful for computation of the codimension growth.
For example, the last result yields an  asymptotic.
\begin{Corr}[\cite{Pe97msb}]
$c_n({\bf{N}}_c{\bf{A}})\approx c^{n-c-1}n^{c}$, as $n\to \infty$.
\end{Corr}

Let $f(z)$  be an entire function of complex variable, denote ${\rm M}_f(r):=\max_{|z|=r}|f(z)|$.
Observe that in case of complexity functions, we have ${\rm M}_f(r)=f(r), r\in {\mathbb R}^+$,
since all coefficients are nonnegative.

By Theorem~\eqref{TcodimNN}, ${\mathcal C}({\bf{N}}_{s_q}\cdots{\bf{N}}_{s_1},z)$
is "almost" $q-1$ iterations of $\exp(*)$ applied to
$\beta_1(z)=\sum_{m=1}^{s_1}{z^m}/m$, thus one has
something like $\exp^{(q-1)}(z^{s_1} /s_1)$.
More, precisely, we derive the following asymptotic.
\begin{Th}[\cite{Pe97msb}]\label{TNNcodimfunc}
Consider the variety of polynilpotent Lie algebras
${\bf{V}}:={\bf{N}}_{s_q}\cdots{\bf{N}}_{s_1}$, $q\ge 2$,
and its complexity function $f(z):={\mathcal C}({\bf{V}},z)$. Then
  $$
  \lim_{r\to\infty}
  \frac{\ln^{(q-1)}{\mathrm M}_f(r)}{r^{s_1}}=\frac{s_2}{s_1}.
  $$
\end{Th}
Next, we establish a relation between the growth of fast growing entire functions and
an asymptotic of their coefficients. This fact helped us to match the upper bounds in the next result.
But a connection between the lower bounds require more direct estimates.
\begin{Th}[\cite{Pe97msb}]\label{TNN2}
Consider the variety of polynilpotent Lie algebras ${\bf{V}}:={\bf{N}}_{s_q}\cdots{\bf{N}}_{s_1}$, $q\ge 2$.
Then there exists an infinitesimal such that
  $$\displaystyle
  c_n({\bf{V}})=
  \left\{
  \begin{array}{ll}
  \displaystyle
  (n!)^{\frac{\scriptstyle s_1-1}{\scriptstyle s_1}}
\big(s_2+o(1)\big)^{n/s_1}\!,\quad
    &q=2;\\[10pt]
  \displaystyle
\frac{n!}{(\ln^{(q-2)}n)^{n/s_1}}\Bigl(\frac{s_2+o(1)}{s_1}
      \Bigr)^{n/s_1}\!,
    &q=3,4,\dots;\quad
  \end{array}
  \right.\quad n\to\infty. $$
\end{Th}
In the following two cases we obtain somewhat more precise asymptotic, but  they look rather complicated.
\begin{Th}[\cite{Pe99JMSciSch_exp}] Fix $d\in \N$. Then
\begin{align*}
c_n({\bf{A}}{\bf{N}}_d)
   &\approx \mu\left(\strut n!\right)^{1-1/d}
     \exp\bigg(\sum_{k=1}^{d-1}\lambda_kn^{1-k/d}\bigg) n{\rule{0pt}{1em}}^{\frac{3-d}{2d}}\!,
  \quad n\to\infty, \quad \text{where}\\
\lambda_k&:=\frac{\left(\frac kd+1\right)\cdots
    \left(\frac kd+k-1\right)}{k!(d-k)},\qquad k=1,\dots,d-1;\\
\mu&:=\exp\bigg(-\frac 1d \sum_{k=2}^d\frac 1k \bigg)
    (2\pi){\rule{0pt}{1em}}^{\frac{1-d}{2d}}d^{-1/2}.
\end{align*}
\end{Th}

It is well known that $c_n({\bf{A}}^2)=n-1\approx n$, this coincides with our asymptotic. The particular cases are.

\begin{align*}
\begin{split}
\displaystyle c_n({\bf{A}}{\bf{N}}_2)
   &\approx
   \sqrt{n!}\,
   \frac{ \exp\left(\sqrt{n}\right)\sqrt[4]{n}}
        { \sqrt[4]{8\pi e}},\\
\displaystyle c_n({\bf{A}}{\bf{N}}_3)
   &\approx
   \left(n!\strut\right)^{2/3}\,
   \frac{
          \exp\left( {\frac 12} n^{2/3}+{\frac 56} n^{1/3}\right) }
        { \sqrt{3} \sqrt[3]{2\pi} e^{5/18} } ,
\end{split}\qquad n\to\infty.
\end{align*}

\begin{Th}[\cite{Pe99JMSciSch_exp}]
Consider the variety ${\bf{A}}^3$ of solvable Lie algebras of length 3 and its codimension growth sequence
$c_n:=c_n({\bf{A}}^3)$. Then
$$ c_n=\frac{n!}{(\ln n)^n}\exp
       \left(\frac n{\ln n}
          \left(
          2\ln\ln n+1+\frac{2(\ln\ln n)^2-2\ln\ln n-1}{\ln n}+
                                  o\left(\frac1 {\ln n}\right)
          \right)
       \right),\quad n\to\infty.
   $$
\end{Th}


\begin{thebibliography}{99}
\bibitem{Ba}
   Bahturin Yu. A., Identical relations in Lie algebras.
   VNU Science Press, Utrecht, 1987.
\bibitem{BMPZ}
   Bahturin~Yu.A.,  Mikhalev~A.A., Petrogradsky~V.M. and Zaicev~M.V.,
   Infinite dimensional Lie superalgebras.
   de Gruyter Exp. Math. {\bf 7}. de Gruyter, Berlin, 1992.
\bibitem{Ber82}
   Berele, A.,
   Homogeneous polynomial identities.
   {\it Israel J. Math}. {\bf 42} (1982), no. 3, 258--272.
\bibitem{Drensky}
  Drensky V., Free algebras and PI-algebras. Graduate course in algebra.
  Springer-Verlag Singapore, Singapore, 2000.
\bibitem{Drensky87}
  Drensky V.,
  Polynomial identities for the Jordan algebra of a symmetric bilinear form,
  {\it J. Algebra} {\bf 108} (1987), 66--87.
\bibitem{Egor}
    Egorychev G.P.,
    {Integral representation and the computation of combinatorial sums},
    Transl. Math. Monogr. vol.~59, Amer. Math. Soc.,
    Providence, RI, 1984.
\bibitem{GiaMiZai08}
  Giambruno, A., Mishchenko, S., Zaicev, M.,
  Codimensions of algebras and growth functions.
  {\it Adv. Math.} {\bf 217} (2008), no. 3, 1027--1052.
\bibitem{GiaZa99}
  Giambruno A., Zaicev M.,
  Exponential codimension growth of PI algebras: an exact estimate.
  {\it Adv. Math.} {\bf 142} (1999), no. 2, 221--243.
\bibitem{GiaZa13}
  Giambruno A., Zaicev M.,
  Non-integrality of the PI-exponent of special Lie algebras.
  {\it Adv. in Appl. Math.} {\bf 51} (2013), no. 5, 619--634.
\bibitem{GiaZai}
   Giambruno A., Zaicev M.
   Polynomial identities and asymptotic methods.
   Mathematical Surveys and Monographs 122. Providence, RI: American Mathematical Society (AMS) (2005).
\bibitem{GiaZe11}
   Giambruno, A., Zelmanov, E.;
   On growth of codimensions of Jordan algebras. in:
   Groups, algebras and applications.
   Providence, RI: AMS. Contemporary Mathematics 537, 205--210 (2011).
\bibitem{Grishkov88}
   Grishkov~A.N.,
   On growth of varieties of Lie algebras,
   {\it Mat. Zametki}, {\bf 44}, No.~1, (1988),  51--54.
   Engl. transl.,
   {\it Math. Notes}, {\bf 44}, (1988), No.~1--2, 515--517.
\bibitem{KlPe05}
   Klementyev\,S.G. and Petrogradsky\,V.M.
   Growth of solvable Lie superalgebras.
   {\it Comm. Algebra.}, {\bf 33} (2005), no.\,3, 865--895.
\bibitem{KlPe05usp}
   Klementyev\,S.G. and Petrogradsky\,V.M.
   On growth of almost solvable Lie algebras.
   {\it Uspekhi Mat. Nauk}, {\bf 60} (2005), no.\,5, 165--166.
   translation in
   {\it Russian Math. Surveys}, {\bf 60} (2005), no.\,5, 970--972.
\bibitem{KlPe09contempmath}
   Klementyev\,S.G. and Petrogradsky\,V.M.
   On growth of almost polynilpotent Lie algebras,
   in {\it Groups, Rings and Group Rings}., ed. A.Giambruno, C.~Polcino~Milies, S.K.~Sehgal,
   Contemp. Math. {\bf 499}, AMS, RI, 2009, 173--180.
   %%% 270 pp., Softcover, ISBN-10: 0-8218-4771-6, ISBN-13: 978-0-8218-4771-8,
\bibitem{KraLen}
   Krause~G.R. and Lenagan~T.H.,
   {Growth of algebras and Gelfand-Kirillov dimension},
   Pitman, Boston, 1985.
\bibitem{KouTet67}
    {Kourovskaya tetrad,  Unsolved problems in group theory},
    Nauka, Novosibirsk, 1967.
\bibitem{Lat63}
  Latyshev V.N.,
  {Two remarks on $PI$-algebras},
  {\it Sibirsk. Mat. \v Z}. {\bf 4} (1963) 1120--1121.
\bibitem{Laz54}
   Lazard M.,
   {\it Sur les groupes nilpotents et les anneaux de Lie},
   {Ann. Sci. \'Ecole Norm. Sup.}
   {\bf 71}, (1954), 101--190.
\bibitem{Licht84}
   Lichtman A.I.,
   {Growth in enveloping algebras},
   {\it Israel J. Math.} {\bf 47}, No.~4, (1984), 297--304.
\bibitem{Pe95}
   Petrogradsky~V.M.,
   On types of overexponential growth of identities in Lie PI-algebras. (Russian)
   {\it Fundam. Prikl. Mat.}~{\bf 1}, no.\,4, (1995), 989--1007.
\bibitem{Pe96}
   Petrogradsky V.M.,
   {Intermediate growth in Lie algebras and their enveloping algebras},
   {\it J.~Algebra} {\bf 179}, (1996), 459--482.
\bibitem{Pe97msb}
   Petrogradsky V.M.,
   Growth of polynilpotent varieties of Lie algebras and rapidly growing entire functions.
   {\it Mat. Sb.}, {\bf 188} (1997), no.\,6, 119--138;
   translation in
   {\it Russian Acad. Sci. Sb. Math}., {\bf 188} (1997), no.\,6, 913--931.
\bibitem{Pe99isr}
   Petrogradsky V.M.,
   Exponential generating functions and complexity of Lie varieties.
   {\it Israel J. Math.} {\bf 113} (1999), 323--339.
\bibitem{Pe99int}
   Petrogradsky  V.M.,
   Growth of finitely generated polynilpotent Lie algebras and groups, generalized partitions, and functions analytic in the unit circle,
   {\it Internat. J. Algebra Comput.}, {\bf 9} (1999), no 2, 179--212.
\bibitem{Pe99JMSciSch_exp}
   Petrogradsky V.M.,
   Exponential Schreier's formula for free Lie algebras and its applications.
   {\it Algebra, {\bf 11}. J. Math. Sci. (New York)},
   {\bf 93} (1999), no.\,6, 939--950.
\bibitem{Pe00archiv}
   Petrogradsky V.M.,
   Schreier's formula for free Lie algebras.
   {\it Arch. Math. (Basel)}, {\bf 75} (2000), no.\,1, 16--28.
\bibitem{Pe00d}
   Petrogradsky V.M.,
  On growth of Lie algebras, generalized partitions, and analytic functions.
  (Formal power series and algebraic combinatorics, Vienna, 1997).
  {\it Discrete Math.}, {\bf 217} (2000), no.\,1--3, 337--351.
\bibitem{Pe02}
   Petrogradsky V.M.,
   On generating functions for subalgebras of free Lie superalgebras.
   {\it Discrete Math.},  {\bf 246} (2002), no.\,1--3,  269--284.
\bibitem{Pe05}
   Petrogradsky~V.M.,
   Enumeration of algebras close to absolutely free algebras and binary trees.
   {\it J.Algebra}, {\bf 290}, (2005), no.\,2, 337--371.
\bibitem{Pe21}
   Petrogradsky V.,
   Scale for codimension growth of Poisson PI-algebras, {\it Israel J. Math.}, to appear, arXiv:2107.02424.
\bibitem{Rats14}
   Ratseev, S. M.
   Correlation of Poisson algebras and Lie algebras in the language of identities.
   {\it Math. Notes} {\bf 96} (2014), no. 3-4, 538--547;
   Translation of {\it Mat. Zametki} {\bf 96} (2014), no. 4, 567--577.
\bibitem{Razmyslov}
   Razmyslov~Yu.P., Identities of algebras and their representations,
   AMS, Providence, RI, 1994.
\bibitem{Reg72}
  Regev~A., Existence of identities in $A\otimes B$.
  {\it Israel J. Math}. {\bf 11} (1972), 131--152.
\bibitem{Shmelkin}
    Shmel'kin A.L.,
    {\it Free polynilpotent groups},
    {Izv. Akad. Nauk SSSR Ser. Mat.} {\bf 28}, No.~1, (1964), 91--122;
    English transl.,
    {Amer. Math. Soc. Transl. ser.}~2, vol.~55,
    Amer. Math. Soc., Providence, RI, 1966, 270--304.
\bibitem{Ufn}
   Ufnarovskiy V.A.,
   {Combinatorial and asymptotic methods in algebra},
   Itogi Nauki i Tekhniki, Sovrem. Probl. Mat. Fund. Naprav.
   vol.~{57}, Moscow, 1989;
   Engl. transl.,
   Encyclopaedia Math. Sci., vol.~{57},
   \newcounter{r6}\setcounter{r6}{6}  Algebra~\Roman{r6},
   Springer, Berlin, 1995.
\bibitem{Vol}
  Volichenko I.~B.,
   {\it On variety $\AA\NN_2$ over field of zero characteristic},
   {Dokl. Akad. Nauk Belarusi}
   \newcounter{r25}\setcounter{r25}{25}\Roman{r25}, No.~12, (1981),
   1063--1066.
\bibitem{Zai02}
   Zaitsev, M.V.,
   Integrality of exponents of growth of identities of finite-dimensional Lie algebras.
   {\it Izv. Ross. Akad. Nauk Ser. Mat.} {\bf 66} (2002), no. 3, 23--48; translation in
   {\it Izv. Math.} {\bf 66} (2002), no. 3, 463--487.
\bibitem{ZaiMi99}
    Zaicev M., Mishchenko S.P.,
    An example of a variety of Lie algebras with a fractional exponent.
    {\it J. Math. Sci.} {\bf 93} (1999), 977--982.
\end{thebibliography}
\end{document}